\documentclass[a4paper,10pt,reqno]{amsart}
\usepackage{amsmath,amsfonts,amssymb,amsthm}
\usepackage{graphics}
\usepackage{epsfig}
\usepackage[margin=2.3cm]{geometry}
\usepackage{esint} 
\usepackage{color}

\newtheorem{theorem}{Theorem}

\newtheorem{lemma}{Lemma}

\providecommand{\ah}{a_{\mathrm{hom}}}

\newcommand{\en}[1]{\left< #1 \right>}
\newcommand{\p}[1]{\left(#1\right)}

\newcommand{\R}{\mathbb{R}}

\newcommand{\ignore}[1]{}

\newcommand{\abs}[1]{\lvert#1\rvert}
\newcommand{\norm}[1]{\lVert#1\rVert}

\author{Peter Bella}
\address{Max Planck Institute for Mathematics in the Sciences, Inselstrasse 22, Leipzig (Germany)}
\email{bella@mis.mpg.de}

\author{Benjamin Fehrman$^\dagger$}\thanks{$^\dagger$The second author is supported by the National Science Foundation Mathematical Sciences Postdoctoral Research Fellowship under Grant Number 1502731.}
\address{Max Planck Institute for Mathematics in the Sciences, Inselstrasse 22, Leipzig (Germany)}
\email{fehrman@mis.mpg.de}

\author{Felix Otto}
\address{Max Planck Institute for Mathematics in the Sciences, Inselstrasse 22, Leipzig (Germany)}
\email{otto@mis.mpg.de}

\title{A Liouville theorem for elliptic systems with degenerate ergodic coefficients}

\begin{document}

\date{May 2, 2016}
\maketitle

\begin{abstract}
  We study the behavior of second-order degenerate elliptic systems in divergence form with random coefficients which are stationary and ergodic.  Assuming moment bounds like Chiarini and Deuschel [Arxiv preprint 1410.4483, 2014] on the coefficient field $a$ and its inverse, we prove an intrinsic large-scale $C^{1,\alpha}$-regularity estimate for $a$-harmonic functions and obtain a first-order Liouville theorem for subquadratic $a$-harmonic functions. 
\end{abstract}

\begin{center}
\begin{Large}
\end{Large}
\end{center}


\section{Introduction and the main results}

We study the behavior of second order non-uniformly elliptic equations, and more generally systems of equations, with random coefficients.  The random coefficient fields $a$ are {\it stationary}, meaning that the joint probability distribution of $a$ and $a(\cdot + x)$ are the same, and {\it ergodic}, meaning that every translation invariant random variable is almost surely constant.  Furthermore, as in the framework of Chiarini and Deuschel~\cite{ChiariniDeuschel2014}, rather than assuming the field is uniformly elliptic we assume only moment bounds from above and below.

More precisely, if $\en{\cdot}$ denotes the expectation with respect to the probability measure on the space of coefficient fields, which will be denoted $\Omega$, we define the scalar random variables $0<\lambda,\mu<\infty$ via
\begin{equation}\nonumber
\lambda:=\inf_{\xi\in\mathbb{R}^d}\frac{\xi\cdot a\xi}{|\xi|^2}\quad\mbox{and}\quad
\mu:=\sup_{\xi\in\mathbb{R}^d}\frac{|a\xi|^2}{\xi\cdot a\xi},
\end{equation}
where in the symmetric case $\lambda^{-1}=|a^{-1}|$ and $\mu=|a|$ are the spectral norms of $a$ and its inverse.  Our assumption is that
\begin{equation}\label{apq}
 \en{ \mu^p }^{\frac{1}{p}} + \en{ \lambda^{-q}}^{\frac{1}{q}} =: K < \infty \quad \textrm{where} \quad\frac{1}{p} + \frac{1}{q} < \frac{2}{d}.
\end{equation}
Here $d \ge 2$ denotes the dimension. Notice that (\ref{apq}) coincides with the integrability condition considered in~\cite{ChiariniDeuschel2014}. 

We require this condition most essentially in two steps of the proof.  First, the strict inequality appearing in (\ref{apq}) is used in the foremost stochastic element of the argument, and guarantees the compactness of a certain Sobolev embedding which is used in Lemma \ref{lm2} to establish the sublinearity of the large-scale averages of the corrector and flux correction defined in (\ref{correq}), (\ref{sigmaeq1}) and \eqref{sigmaeq2}. The deterministic elements of the argument, namely the Caccioppoli inequality of Lemma \ref{lmcac} and the large-scale $C^{1,\alpha}$-regularity of Theorem~\ref{prop1}, require only $1/p + 1/q \leq 2/d$.

The primary result of this paper is a {\it first-order} Liouville theorem for degenerate coefficient fields satisfying \eqref{apq}. By which we mean that for $\en{\cdot}$-a.e environment $a$, every subquadratically-growing $a$-harmonic function $u$ on the whole space belongs to the~$(d+1)$-dimensional space of $a$-affine functions. Namely, the space spanned by functions of the form  $c + \xi\cdot x + \phi_\xi$ where, for every $\xi\in\mathbb{R}^d$, the corrector $\phi_\xi$ denotes the solution
\begin{equation}\label{correq}
 -\nabla \cdot a (\xi + \nabla \phi_\xi) = 0,
\end{equation}
whose gradient $\nabla \phi_\xi$ is {\it stationary}, by which we understand $\nabla \phi_i(a;x+z) = \nabla \phi_i(a(\cdot+z),x)$ for any shift vector $z \in \R^d$, has vanishing average $\en{\nabla \phi_\xi} = 0$, and has finite second moment $\en{|\nabla \phi_\xi|^2}$.  The following theorem summarizes the result.

\begin{theorem}\label{liouville}
 Let $\en{\cdot}$ be stationary and ergodic, and assume it satisfies~\eqref{apq}.  Then $\en{\cdot}$-a.e. coefficient field $a \in \Omega$ has the following Liouville property:
 if $u$ is an $a$-harmonic function in the whole space, i.e., it solves $- \nabla \cdot a \nabla u = 0$ in $\R^d$, and in addition $u$ is subquadratic in the sense that, for some $\alpha<1$,
 \begin{equation}\nonumber
  \lim_{R \to \infty} R^{-(1+\alpha)} \biggl(\fint_{B_R} |u|^\frac{2p}{p-1}\biggr)^{\frac{p-1}{2p}} = 0,
 \end{equation}
 then there necessarily exists $c\in\mathbb{R}$ and $\xi\in\mathbb{R}^d$ for which $u(x) = c + \xi\cdot x + \phi_\xi(x)$. 
\end{theorem}

Theorem~\ref{liouville} amounts to an extension of the first-order Liouville property for uniformly elliptic coefficient fields obtained by Gloria, Neukamm, and the third author~\cite[Corollary~1]{GNO4} to the case of degenerate environments satisfying (\ref{apq}).  For this, like in \cite{GNO4}, we will also need, for each $i\in\{1,\ldots,d\}$, the skew-symmetric matrix field
$\sigma_i$, which can be viewed as a vector potential for the harmonic coordinates and provides the correction of the $i$th component of the flux
$$q_i := a(\nabla \phi_i + e_i),$$
by satisfying the equation
\begin{equation}\label{sigmaeq1}
 q_i - \en{q_i} =: \nabla \cdot \sigma_i.
\end{equation}
Here the divergence of a tensor field is defined as $$\p{ \nabla \cdot \sigma_i }_j := \sum_{k=1}^d \partial_k \sigma_{ijk}.$$
The linearity of the equation allows for the consideration of only the fluxes corresponding to the canonical basis $\{e_i\}_{i\in\{1,\ldots,d\}}$, where $$\en{q_i} = \en{a(\nabla\phi_i+e_i)}=:\ah e_i$$ defines the constant, possibly non-symmetric, homogenized coefficients $\ah$.  The uniform ellipticity of $\ah$ was established by \cite[Proposition~4.1]{ChiariniDeuschel2014}.

In the setting of uniformly elliptic random coefficient fields the vector corrector $\sigma$ was introduced and constructed in~\cite[Lemma~1]{GNO4}. Since the definition of $\sigma_i$ is under-determined, taking motivation from the analogous periodic framework, they made the specific choice of gauge
\begin{equation}\label{sigmaeq2}
 -\Delta \sigma_{ijk} = \partial_j q_{ik} - \partial_k q_{ij}.
\end{equation}
A principle difference between the degenerate and uniformly elliptic cases is that, in the latter the fluxes $q_i$ belong to $L^2(\Omega)$, and therefore so too do the gradients $\nabla \sigma_{ijk} \in L^2(\Omega)$, whereas in the degenerate setting this is no longer true. In fact, H\"older's inequality and the moment bound (\ref{apq}) imply only that $q_i \in L^\frac{2p}{p+1}(\Omega)$. For this reason, it is necessary to replace the $L^2$-theory for the construction of $\sigma$ used in \cite[Lemma~1]{GNO4} with an approximation argument and a Calder\'on-Zygmund estimate to first construct the stationary, mean zero gradients of $\sigma$ in $L^\frac{2p}{p+1}(\Omega)$, and thereby define the generally non-stationary flux corrections uniquely up to an additive random constant.

The properties of the correctors and flux corrections are summarized by the following lemma, where we remark that the construction of the scalar corrector $\phi$ in the degenerate ergodic setting under the weaker assumptions $p=q=1$ has been carried out in~\cite[Section 4]{ChiariniDeuschel2014}. Loosely speaking, and for the construction of both the corrector and the flux corrections, the definitions (\ref{correq}), (\ref{sigmaeq1}) and (\ref{sigmaeq2}) are lifted to the probability space in order to construct their gradients as stationary, mean-zero, finite-energy random fields.

\begin{lemma}\label{lm1}
 Let $\en{\cdot}$ be stationary and ergodic, and let~\eqref{apq} be satisfied. Then there exist $C=C(d)>0$ and two random tensor fields $\{ \phi_i \}_{i=1,\ldots,d}$ and $\{\sigma_{ijk}\}_{i,j,k=1,\ldots,d}$ with the following properties: 
 The gradient fields are stationary, have bounded moments, and are of vanishing expectation:
 \begin{equation}\label{c1}
   \sum_{i=1}^d \en{ \nabla\phi_i\cdot a\nabla \phi_i } + \sum_{i=1}^d \en{ |\nabla \phi_i|^{\frac{2q}{q+1}} }^\frac{q+1}{2q} + \sum_{i,j,k=1}^d \en{ |\nabla \sigma_{ijk} |^\frac{2p}{p+1} }^{\frac{p+1}{2p}} \le CK, \quad \en{\nabla \phi_i} = \en{\nabla \sigma_{ijk} } = 0.
 \end{equation} 
Moreover, the field $\sigma$ is skew symmetric in its last two indices, that is
\begin{equation}\nonumber
 \sigma_{ijk} = -\sigma_{ikj}.
\end{equation}
Furthermore, for $\en{\cdot}$-a.e. $a$ we have
\begin{equation}\nonumber
 q_i = a(\nabla \phi_i + e_i) = \ah e_i + \nabla \cdot \sigma_i.
\end{equation}
Finally, the homogenized coefficient field $\ah$ 
is uniformly elliptic in the sense that, for each $\xi\in\mathbb{R}^d$, $$\frac{1}{K}\abs{\xi}^2\leq \xi\cdot \ah \xi\;\;\textrm{and}\;\;\abs {\ah \xi}\leq K\abs{\xi}.$$
\end{lemma}

Furthermore, owing to the fact that the gradients of the corrector $\phi$ and flux-correction $\sigma$ have zero average, the ergodicity of the coefficient field guarantees by standard arguments that their large-scale averages are sublinear in the sense of the following lemma.

\begin{lemma}\label{lm2}  Let $\en{\cdot}$ be stationary and ergodic, and let~\eqref{apq} be satisfied.  Then, the large-scale averages of the random tensor fields $\{\phi_i\}_{i=1,\ldots,d}$ and $\{\sigma_{ijk}\}_{i,j,k=1,\ldots,d}$ are sublinear in the sense that, for $\en{\cdot}-a.e.$ a,  $$\begin{aligned}
  \lim_{R\rightarrow \infty}\frac{1}{R} \biggl( \fint_{B_R} |\phi - \fint_{B_R} \phi|^{\frac{2p}{p-1}} \biggr)^{\frac{p-1}{2p}} & =0,
  \\
 \lim_{R\rightarrow\infty} \frac{1}{R} \biggl( \fint_{B_R} |\sigma - \fint_{B_R} \sigma|^{\frac{2q}{q-1}} \biggr)^{\frac{q-1}{2q}} &=0.
  \end{aligned}$$\end{lemma} \noindent  We remark that an alternate construction of the flux correction is presented in the appendix, and an ingredient of this argument requires a small modification of Lemma \ref{lm1}.  Indeed, the proof of sublinearity follows from the integrability of the gradient fields $\nabla\phi$ and $\nabla\sigma$ and does not use any properties of the underlying equation.
  
The large-scale $C^{1,\alpha}$-regularity first obtained in \cite{GNO4} asserts that whenever $u$ is an $a$-harmonic function, then its deviation from the space of $a$-harmonic affine functions, as defined for each $r>0$ by the excess
$$\mathrm{Exc}(r)=\inf_{\xi\in\mathbb{R}^d}\fint_{B_r}(\nabla u-(\xi+\phi_\xi))\cdot a(\nabla u-(\xi+\phi_\xi)),$$
decays for all sufficiently large radii and any $\alpha\in(0,1)$ as a power law in $(r/R)^{2\alpha}$. The proof is purely deterministic and is based on estimating the homogenization error determined by an $a$-harmonic function $u$ and an $\ah$-harmonic function $v$, as defined by
$$u-(v+\phi_i\partial_iv),$$
for $\phi_i$ the first-order corrector defined in (\ref{correq}) corresponding to the $i$th standard basis vector.  An essential observation of \cite{GNO4} was that the homogenization error satisfied a divergence-form equation with right-hand side in divergence-form.  We use this fact to estimate its energy in the intrinsic $L^2(a)$-norm, where the regularity of the $\ah$-harmonic function $v$ plays an essential role, and to ultimately prove the excess decay and large-scale $C^{1,\alpha}$-regularity.

In this setting, the construction of the appropriate $\ah$-harmonic function $v$ differs considerably from the uniformly elliptic case.  To estimate the homogenization error on the ball $B_R$, the idea is to exploit the best integrability of the coefficient field by separating
\begin{equation}\label{cases_intro}
\mbox{the ``Dirichlet case''}\;q\ge p\;\;\mbox{and the}\;\mbox{``Neumann case''}\;p\ge q,
\end{equation}
where in the Dirichlet case, we define $v$ via the boundary condition
\begin{equation}\nonumber
v=u_\epsilon\;\mbox{on}\;\partial B_R,
\end{equation}
and, in the Neumann case, we impose
\begin{equation}\nonumber
\nu\cdot \ah\nabla v=(\nu\cdot a\nabla u)_\epsilon\;\mbox{on}\;\partial B_R,
\end{equation}
where the subscript $\epsilon$ denotes a smoothing by convolution on the boundary of the ball.  Then, like in~\cite{GNO4}, the energy of the corresponding homogenization error is controlled by introducing a cutoff $\eta$ vanishing near the boundary and estimating the intrinsic energy of the quantity
$$u-(v+\eta\phi_i\partial_iv),$$
where it will be necessary to use the aforementioned divergence-form equation satisfied by the homogenization error, as modified by the introduction of the cutoff, and to control the subsequent boundary terms arising from the case (\ref{cases_intro}) and the vanishing of $\eta$ on $\partial B_R$.  The result is summarized by the following deterministic theorem, where the constants $C_0$ and $C_1$ depend upon $K$ from (\ref{apq}) through the ellipticity of $\ah$ appearing in Lemma \ref{lm1}.

\begin{theorem}\label{prop1}  Let the H\"older exponent $\alpha \in (0,1)$ and $\Lambda > 0 $ be given. Then there exist constants $C_0,C_1=C_0,C_1(d,\alpha,K,\Lambda)$ with the following property:
 
 If $r < R$ are two radii such that for any $\rho \in [r,R]$ we have 
 \begin{equation}\label{abound}
  \biggl( \fint_{B_\rho} \mu^p \biggr)^{\frac{1}{p}} + \biggl( \fint_{B_\rho} \lambda^{-q} \biggr)^{\frac{1}{q}} \le \Lambda,   
 \end{equation}
 with the exponents $p$ and $q$ satisfying
 \begin{equation}\label{pq}
  \frac{1}{p} + \frac{1}{q} \le \frac{2}{d},
 \end{equation}
 and 
 \begin{equation}\label{corrsmall}
  \begin{aligned}
  \frac{1}{\rho} \biggl( \fint_{B_\rho} |\phi - \fint_{B_\rho} \phi|^{\frac{2p}{p-1}} \biggr)^{\frac{p-1}{2p}} &\le \frac{1}{C_0},
  \\
  \frac{1}{\rho} \biggl( \fint_{B_\rho} |\sigma - \fint_{B_\rho} \sigma|^{\frac{2q}{q-1}} \biggr)^{\frac{q-1}{2q}} &\le \frac{1}{C_0},
  \end{aligned}
 \end{equation}
 then any $a$-harmonic function $u$ in $B_R$, i.e., weak solution of $- \nabla \cdot a \nabla u = 0$ in $B_R$, satisfies
 \begin{equation}\nonumber
  \mathrm{Exc}(r) \le C_1 \biggl( \frac{r}{R} \biggr)^{2\alpha} \mathrm{Exc}(R),  
 \end{equation}
 where the excess 
 \begin{equation}\nonumber
\mathrm{Exc}(\rho) := \inf_{\xi \in \R^d} \fint_{B_\rho} (\nabla u-(\xi+\phi_\xi))\cdot a (\nabla u - (\xi + \nabla \phi_\xi))
 \end{equation}
 measures in the $L^2(a)$-sense deviations of $u$ from the set of $a$-affine functions. 
\end{theorem}

We remark that the assumptions of Theorem \ref{prop1} will be satisfied for $\en{\cdot}$-a.e. environment, provided the radius $r$ is chosen sufficiently large.  Indeed, for any $\alpha\in(0,1)$ and any $C_0>0$, the ergodic theorem asserts that there exists for $\en{\cdot}$-a.e. environment a random radius $r_1=r_1(a)$ such that (\ref{abound}) is achieved for $\Lambda = 2(\en{ \mu^p } + \en{ \lambda^{-q} })$ whenever $r\geq r_1$ and Lemma \ref{lm2} guarantees the existence of $r_2=r_2(a)$ such that (\ref{corrsmall}) is satisfied for every $r\geq r_2$.

Finally, a version of the the Caccioppoli inequality adapted to the degenerate setting will be used in the proofs of Theorems \ref{liouville} and \ref{prop1}.  In the uniformly elliptic case, the statements may be used to bound the $L^2$-norm of the gradient of an $a$-harmonic function on ball by the $L^2$-norm of the function itself on a somewhat larger ball.  A straightforward modification yields the analogous statement for elliptic systems with non-symmetric degenerate coefficients.

\begin{lemma}\label{lmcac} Suppose that $u$ is an $a$-harmonic function on $B_R$, and that for some exponents $p \in (1,\infty)$, $q \in [1,\infty)$ we have
\begin{equation}\label{cac1}
\biggl( \fint_{B_R} \mu^p \biggr)^\frac{1}{p} + \biggl( \fint_{B_R} \lambda^{-q} \biggr)^\frac{1}{q} \le \Lambda.
\end{equation}
Then there exists $C_1, C_2=C_1, C_2(d)>0$ such that for any $0<\rho<\frac{R}{2}$ and any $c\in\mathbb{R}$,
\begin{equation}\label{Cac}
\biggl(\fint_{B_{R-\rho}}\abs{\nabla u}^{\frac{2q}{q+1}}\biggr)^{\frac{q+1}{q}}\leq C_1\Lambda \fint_{B_{R-\rho}}\nabla u\cdot a\nabla u \leq C_2\frac{\Lambda^2}{\rho^2}\biggl(\fint_{B_R\setminus B_{R-\rho}}\abs{u-c}^{\frac{2p}{p-1}}\biggr)^{\frac{p-1}{p}}.\end{equation}
\end{lemma}

In the uniformly elliptic framework, the Caccioppoli inequality~\eqref{Cac} can be viewed as a version of a reversed Poincar\'e inequality, meaning that we gain one derivative in the estimate at the expense of increasing the radius of the ball. With the assumption of uniform ellipticity replaced by a weaker moment bound condition on $a$ from below and above, one has to replace the integrability exponents in~\eqref{Cac} on both sides. Hence in this case one trades a derivative for a possible loss in the integrability. While in Lemma~\ref{lmcac} we did not assume condition~\eqref{pq}, which appeared in Theorem~\ref{prop1}, it has a direct relation to~\eqref{Cac}. Indeed, if one uses Sobolev embedding on the right-hand side of~\eqref{Cac} to trade one derivative for better integrability, it is exactly condition~\eqref{pq} which ensures that in the end we get the same exponent as the one we started with on the left-hand side of~\eqref{Cac}. In the case of a condition on $p$ and $q$ with strict inequality~\eqref{apq}, the above combination of Caccioppoli and Sobolev inequalities gives a gain in the integrability -- a fact that allowed Chiarini and Deuschel~\cite{ChiariniDeuschel2014} (see also~\cite{ChiariniDeuschel2015}), in the case of a scalar equation, to perform a Moser iteration. The condition~\eqref{apq} first appeared in the paper by Andres, Deuschel, and Slowik~\cite{AndresDeuschelSlowik2015} (see also~\cite{AndresDeuschelSlowik2015Heat}), and was recently generalized to study invariance principles for environments with time-dependent coefficients~\cite{AndresChiariniDeuschelSlowik2016,DeuschelSlowik2016}.

First-order Liouville statements for $a$-harmonic functions are a compact way to express regularity on large scales.
In fact, an easy post-processing of the excess decay in Theorem~\ref{prop1} yields a large-scale $C^{1,\alpha}$-estimate for
$a$-harmonic functions, see~\cite[Corollary 2]{GNO4}. We thus speak of a $C^{1,\alpha}$-Liouville property. A further post-processing
yields large-scale $C^{1,\alpha}$-Schauder estimates for the operator $-\nabla\cdot a\nabla$, see for instance~\cite[Theorem 5.20]{giaqmart}. In the case of constant-coefficient (and thus scale-invariant) equations, this relation between $C^{1,\alpha}$-Liouville principles and $C^{2,\alpha}$-Liouville principles on the one hand, namely that sub-cubic harmonic functions must be quadratic harmonic polynomials, and a $C^{1,\alpha}$- or $C^{2,\alpha}$-Schauder  theory on the other hand is classical: An indirect argument by Simon~\cite{Simon97} allows to directly pass from the Liouville property to the corresponding Schauder estimate.

For general non-constant coefficient fields $a$, we call $C^{k,\alpha}$-Liouville property the fact that
the linear space of $a$-harmonic functions that grow at most of the order $|x|^{k+\alpha}$ (say, in an averaged sense as
in Theorem~\ref{liouville}) has the same dimension as in the case of constant-coefficient equations (where it is spanned by all harmonic polynomials
of order at most $k$). Without further structural conditions, it is almost folkloric knowledge that this equality already fails for $k=0$ and any $\alpha>0$ even in the case of uniformly elliptic coefficients (which may even be smooth~\cite[Proposition 21]{FischerOtto}).  The work of Yau~\cite{YauHarmonic}, drawing a connection to curvature of the metric given by $a$, 
popularized the question of determining whether the dimensions are {\it asymptotically} equal for $k\uparrow\infty$, as shown by Colding and Minicozzi~\cite{ColdingMinicozzi} and Li~\cite{Li1997} for uniformly elliptic equations.

In the case of uniformly elliptic {\it periodic} coefficient fields, 
the full hierarchy of Liouville properties was established by Avellaneda and Lin \cite{AvellanedaLinCRAS89}, based on earlier ideas developed by those authors on a large-scale regularity theory in H\"older and $L^p$-spaces~\cite{AvellanedaLinCPAM87,AvellanedaLin87} via a Campanato iteration, which is also used in Theorem~\ref{prop1}. 
Marahrens and the last author~\cite[Corollary 4]{MO} derived a $C^{0,\alpha}$-Liouville property in the case of stationary {\it random} coefficient fields with integrable correlation tails (that is, integrable in a sufficiently strong sense so as to allow for a Logarithmic Sobolev Inequality).
Benjamini, Duminil-Copin, Kozma and Yadin \cite{BenjaminiDuminilKozma+} derived a $C^{0,\alpha}$-Liouville property under the mere assumption
of ergodicity, and that allows for degenerate coefficient fields provided suitable heat-kernel bounds are
available. They also formulated the question of higher-order Liouville properties in the random case \cite[Theorems~4,5]{BenjaminiDuminilKozma+}.
Armstrong and Smart \cite{armstrongsmart2014} adapted the approach of Avellaneda and Lin \cite{AL1} to obtain a large-scale $C^{1,0}$-regularity theory
in the case of uniformly elliptic coefficient fields with a finite-range condition, which was a major step because it required a new quantitative substitute for the compactness argument, and which was later extended by Armstrong and Mourrat~\cite{ArmstrongMourrat}
to very general mixing conditions.  Gloria, Neukamm and the last author \cite{GNO4} derived the $C^{1,\alpha}$-Liouville property under the mere assumption of ergodicity
in the uniformly elliptic case; the main new ingredients being 1) the usage of an {\it intrinsic} excess decay, that is, 
measuring the energy distance to the space of {\it intrinsically} affine functions (i.e., the harmonic coordinates) 
and 2) the construction of the vector potential $\sigma$ (which allows to bring the residuum in the two-scale expansion
into divergence form). Fischer and the last author extended Theorem 2 to the case of an excess of order $k$ under a mild quantification of the sublinear growth of the corrector in~\cite{FischerOtto} to obtain a full hierarchy of Liouville properties, and showed in~\cite{FischerOtto2} that the quantified sublinear growth of the corrector is satisfied under a mild quantification of ergodicity in a certain class of Gaussian environments.  Additionally, there has recently been a lot of activity aimed at understanding the space of harmonic functions on infinite percolation clusters with specified polynomial growth.  Recently, for instance, Sapozhnikov \cite{Artem} proved the finite-dimensionality of these spaces for a large class of percolation models.


Finally, we believe these results are very likely extendable to the discrete case.  Indeed, Deuschel, Nguyen and Slowik \cite{DeuschelNguyenSlowik} have established an invariance principle for random walks in a degenerate environments under similar integrability assumptions on the coefficient field, and the techniques of this paper are expected to be adaptable to their setting.

{\textbf{Organization and Notation.}}  The remainder of the paper presents the proofs of Theorem \ref{liouville}, Lemmas \ref{lm1} and \ref{lm2}, Theorem \ref{prop1} and Lemma \ref{lmcac} in the order of their appearance in the introduction.  In addition, an appendix contains an alternative argument for the construction of the flux corrector $\sigma$ in Lemma \ref{lm1}.  We remark that, in order to simplify the notation, the statements and proofs are written for the non-symmetric {\it scalar} setting.  However, at the cost of increasing some constants, all of the arguments carry through unchanged for non-symmetric systems.  Throughout, $\lesssim$ is used to denote a constant whose dependencies are specified in every case by the statement of the respective lemma or theorem.


\section{The Proof of Theorem \ref{liouville}}

Suppose that $u$ is an $a$-harmonic function on the whole space, that is 
\begin{equation}\nonumber
-\nabla\cdot a\nabla u=0\;\; \textrm{on}\;\;\mathbb{R}^d, 
\end{equation}
which is strictly subquadratic in the sense that, for some $\alpha<1$, 
\begin{equation}\nonumber
\lim_{r\rightarrow\infty}r^{-(1+\alpha)}\biggl(\fint_{B_r}\abs{u}^{\frac{2p}{p-1}}\biggr)^{\frac{p-1}{2p}} = 0. 
\end{equation}
For $\en{\cdot}$-a.e. $a$ it follows from the ergodic theorem and the integrability assumption \eqref{apq} that there exists $r_1=r_1(a)$ such that, for all $r\geq r_1$, 
\begin{equation}\nonumber
\biggl( \fint_{B_r} \mu^p \biggr)^\frac{1}{p} + \biggl( \fint_{B_r} \lambda^{-q} \biggr)^\frac{1}{q} \leq 2({\en{\mu^p}^\frac{1}{p} +\en{\lambda^{-q}}}^\frac{1}{q})=:\Lambda.
\end{equation}
 Let $C_0 = C_0(d,\alpha,K,\Lambda)$ be as in Theorem~\ref{prop1}, and choose $r_2=r_2(a)\geq r_1$ so that, in view of Lemma \ref{lm2}, for all $r\geq r_2$,  \begin{equation}\nonumber
  \frac{1}{r} \biggl( \fint_{B_r} |\phi - \fint_{B_r} \phi|^{\frac{2p}{p-1}} \biggr)^{\frac{p-1}{2p}} \le \frac{1}{C_0}\;\;\textrm{and}\;\; \frac{1}{r} \biggl( \fint_{B_r} |\sigma - \fint_{B_r} \sigma|^{\frac{2q}{q-1}} \biggr)^{\frac{q-1}{2q}} \le \frac{1}{C_0}.
 \end{equation}
 
In order to conclude, observe that Lemma \ref{lmcac} and the definition of excess imply by the choice of $r_1$ that, for each $r\geq r_1$, for $C_1=C_1(d)>0$,
\begin{equation}\nonumber
\mathrm{Exc}(r)\leq\fint_{B_{r}}\nabla u\cdot a\nabla u\le \frac{C_1\Lambda}{r^2}\biggl(\fint_{B_{2r}}\abs{u}^{\frac{2p}{p-1}}\biggr)^{\frac{p-1}{p}}.
\end{equation}
This implies that, in view of Theorem \ref{prop1} and the choice of $r_2$, for every of $r >\rho>r_2$, for $C_2, C_3=C_2, C_3(d,\alpha,K,\Lambda)>0$,
\begin{equation}\nonumber
\mathrm{Exc}(\rho)\leq C_2\left(\frac{\rho}{r}\right)^{2\alpha}\mathrm{Exc}(r)\le C_3 \rho^{2\alpha} r^{-(2+2\alpha)}\biggl(\fint_{B_{2r}}\abs{u}^{\frac{2p}{p-1}}\biggr)^{\frac{p-1}{p}}. 
\end{equation}
Therefore, owing to the choice of $\alpha$, we have, for each $\rho>r_2$,
\begin{equation}\nonumber
\mathrm{Exc}(\rho)\le C_3\rho^{2\alpha} \limsup_{r\rightarrow\infty} \biggl( r^{-(1+\alpha)}\biggl(\fint_{B_{2r}}\abs{u}^{\frac{2p}{p-1}}\biggr)^{\frac{p-1}{2p}}\biggr)^2=0.
\end{equation}
By definition of $\mathrm{Exc}(\rho)$ this implies existence of $\xi \in \R^d$ and $c \in \R$ s.t. $u(x) = c + \xi \cdot x + \phi_\xi(x)$ for a.e. $x \in B_\rho$. Since the values $\xi$ and $c$ are independent of the choice of $\rho$, we obtain the statement of the theorem.

\section{Proof of Lemma~\ref{lm1}}

 The construction of the corrector $\phi$ in the case of degenerate and unbounded, stationary and ergodic coefficients was performed in~\cite[Section 4.1]{ChiariniDeuschel2014}. 

 For the construction of the flux corrector $\sigma$, we combine an approximation argument with a version of the existence result for~$\sigma$ from~\cite{GNO4}. Given $n \in \mathbb{N}$ let us consider the random variable $q_n := I(|q_n| \le n)q$.  Here, $I$ stands for the characteristic (indicator) function. We will prove the existence of a random tensor field $\sigma_n$ which satisfies
 \begin{itemize}
  \item $\nabla \sigma_{n,ijk} \in L^2(\Omega)$ is stationary, $\en{\nabla \sigma_{n,ijk} } = 0$, and $\sigma_n$ is skew-symmetric in its last two indices,
  \item for $\en{\cdot}$-a.e. $a$ we have
  \begin{align}\label{lm1_11}
  \nabla\cdot \sigma_{n,i}=q_{n,i}-\en{q_{n,i}};
  \end{align}
  \item for $\en{\cdot}$-a.e. $a$ we have 
 \begin{align}\label{lm1_1}
 -\Delta \sigma_{n,ijk} = \partial_j q_{n,ik} - \partial_k q_{n,ij}.
 \end{align}
 \end{itemize}
This fact follows from the argument of Gloria, Neukamm, and the third author~\cite[Lemma 1]{GNO4}; for the reader's convenience we outline here an alternative approach which follows the reasoning presented by the third author at the September, 2015 Oberwolfach workshop on stochastic homogenization.

Fix $n > 0$.  The argument first constructs the gradient of the expected approximate flux correction $\sigma_n=\{\sigma_{n,ijk}\}_{i,j,k=1,\ldots,d}$ by considering the single component $\sigma_{n,i}=\{\sigma_{n,ijk}\}$ for each $i\in\{1,\ldots,d\}$ separately.  For this, the equation will be lifted to the probability space, and phrased in terms of the ``horizontal gradient'' with respect to shifts of the coefficient field.

Precisely, for each $i\in\{1,\ldots,d\}$, the horizontal derivative of a random variable $\zeta$ along the $i$th coordinate direction is defined by the infinitesimal generator of the corresponding translation in the probability space, and is given by the limit
$$D_i\zeta(a):=\lim_{h\rightarrow 0}\frac{\zeta(a(\cdot+he_i))-\zeta(a)}{h}.$$
We remark that the operators $D_i$ are closed, and densely defined on $L^2(\Omega)$.  We write $\mathcal{D}(D_i)$ for their respective domains and define the Hilbert space
$$\mathcal{H}^1=\cap_{i=1}^d\mathcal{D}(D_i)\subset L^2(\Omega)\;\textrm{with inner product}\;(f,g)_{\mathcal{H}^1}:=\en{fg}+\sum_{i=1}^d\en{D_if D_ig}.$$
The space $\mathcal{H}^1$ will be used to lift the weak formulations of (\ref{lm1_11}) and (\ref{lm1_1}) to the probability space and ultimately to construct the horizontal gradient of the approximate flux correction.

Henceforth, we fix $i\in\{1,\ldots,d\}$ and $n\geq 0$, and to simplify notation suppress the dependence on both indices in the argument to follow.  Consider the closed subspace of $L^2(\Omega)$ defined by
\begin{equation}\nonumber 
X=\{\;\{S_{ljk}\}_{l,j,k=1,\ldots,d}\in L^2(\Omega;\mathbb{R}^{d^3})\;|\;S_{ljk}+S_{lkj}=0, \partial_mS_{ljk}=\partial_l S_{mjk}\;\textrm{and}\;\en{S_{ijk}}=0\},
\end{equation}
where for a random variable $\zeta$ and for every $i\in\{1,\ldots,d\}$, the notation $\partial_i\zeta$ denotes the distributional derivative of $\zeta$ defined by
$$\en{\partial_i\zeta\chi}=-\en{\zeta D_i\chi}\;\;\textrm{for every}\;\;\chi\in\mathcal{H}^1.$$
We observe that $X$ is a Hilbert space with respect to the standard inner product on $L^2(\Omega;\mathbb{R}^{d^3})$ and that, formally,
we expect the gradient $\{\partial_l\sigma_{ijk}\}_{l,j,k=1,\ldots,d}$ to be an element of $X$ where $S_{ljk}:=\partial_l\sigma_{ijk}.$

Interpreting equation (\ref{lm1_1}) on the space $X$, Riesz' representation theorem and the boundedness of $q$ yield a unique element $\{\overline{S}_{ljk}\}\in X$ satisfying
\begin{equation}\label{lm1_3}
\en{\overline{S}_{jkl} S_{jkl}}=-2\en{q_{k} S_{jjk}}\;\;\textrm{for every}\;\;\{S_{ljk}\}_{l,j,k=1,\ldots,d}\in X,
\end{equation}
where we have employed Einstein's summation convention and $\en{\cdot}$ denotes the standard inner product on $L^2(\Omega)$.

In order to verify (\ref{lm1_1}), it is necessary to prove that, in the sense of distributions, and again employing Einstein's summation convention, for each $j,k\in\{1,\ldots,d\}$,
\begin{equation}\label{lm1_5}-\partial_l\overline{S}_{ljk}=\partial_j q_{k}-\partial_k q_{j}.\end{equation}
As mentioned above, for any skew-symmetric $\{\eta_{jk}\}_{j,k=1,\ldots,d}\in \mathcal{H}^1$ the gradient satisfies $\{D_l\eta_{jk}\}_{l,j,k=1,\ldots,d}\in X$.  Therefore, for an arbitrary such $\{\eta_{jk}\}_{j,k=1,\ldots,d}$, equation (\ref{lm1_3}) implies that 
\begin{equation}\nonumber 
\en{\overline{S}_{ljk}D_l\eta_{jk}}=-2\en{q_{k}D_j\eta_{jk}}=-\en{q_{k}D_j\eta_{jk}}-\en{q_{j}D_k\eta_{kj}}=-\en{q_{k}D_j\eta_{jk}}+\en{q_{j}D_k\eta_{jk}},
\end{equation}  
which, since the skew-symmetric $\{\eta_{jk}\}_{j,k=1,\ldots,d}$ was arbitrary and such functions are dense in $X$, completes the proof of~\eqref{lm1_5}. 

It remains to prove (\ref{lm1_11}) which, when interpreted on the space $X$ turns for each $j\in\{1,\ldots,d\}$ into
\begin{equation}\label{lm1_6} \overline{S}_{kjk}=q_{j}-\en{q_{j}}.\end{equation}
And for this, since $\en{\overline{S}_{ljk}}=0$ for every $l,k,j\in\{1,\ldots,d\}$, the ergodicity implies that it is sufficient to prove that, in the sense of distributions,
\begin{equation}\label{lm1_7}\partial_l\partial_l(\overline{S}_{kjk}-q_{j})=0.\end{equation}
But this follows immediately from the properties of $X$ and (\ref{lm1_5}), which provide the distributional equality
$$\partial_l\partial_l\overline{S}_{kjk}=\partial_l\partial_k\overline{S}_{ljk}=\partial_k\partial_l\overline{S}_{ljk}\stackrel{(\ref{lm1_5})}{=}\partial_k\partial_kq_{j}-\partial_k\partial_jq_{k}=\partial_k\partial_kq_{j}-\partial_j\partial_kq_{k}=\partial_l\partial_lq_{j},$$
where the final inequality is obtained using the fact that $q$ is divergence free.  This completes the argument for (\ref{lm1_7}) and therefore (\ref{lm1_6}) as well.

To conclude, recalling that $i\in\{1,\ldots,d\}$ was fixed throughout, the gradient is defined for each $l,j,k\in\{1,\ldots,d\}$ as $$\partial_l\sigma_{n, ijk}:=\overline{S}_{ljk},$$ which in turn defines each component of the flux correction $\sigma_{n,i}$ and therefore the flux correction $\sigma_n$ itself uniquely up to a random but spatially constant, skew-symmetric vector.  This finishes the proof of existence.

To complete the proof of the lemma, it suffices to prove the uniform in $n$ estimates for the expectation $\en{ |\nabla \sigma_{n,ijk}|^\frac{2p}{p+1} }$.  The result then follows by taking the limit $n \to \infty$. More generally, given two random fields $f \in L^\infty(\Omega;\R^d)$ and $\sigma$, such that $\en {\nabla \sigma} = 0$, $\nabla \sigma \in L^2(\Omega)$ is stationary, $\sigma$ is skew-symmetric, and $\sigma$ and $f$ are related through 
 \begin{equation}\label{p003}
  -\Delta \sigma = - \nabla \cdot f,
 \end{equation}
 it is enough to show a Calder\'on-Zygmund type estimate 
 \begin{equation}\label{p004}
  \en{ |\nabla \sigma|^r } \le C(d,r) \en{ |f|^r },
 \end{equation}
 for general $1 < r < \infty$. 

 For $R,T > 0$ we consider $\sigma_{T,R}$, an approximation of $\sigma$, defined as a unique finite energy solution of
 \begin{equation}\nonumber
  \tfrac{1}{T} \sigma_{T,R} - \Delta \sigma_{T,R} = - \nabla \cdot (\eta_R f),
 \end{equation}
 where $\eta_R$ is a radial cut-off function for $\{ |x| < R \}$ in $\{ |x| < 2R \}$. The additional term $\tfrac{1}{T}  \sigma_{T,R}$ is called massive term, and localizes (up to an exponentially decay) the spatial dependence of the solution on the right-hand side. In the physics community, the above equation is called a \emph{screened Poisson equation}. By the standard Calder\'on-Zygmund theory of singular integral operators, applied to the massive Green's function (in fact its second mixed derivative), we get an estimate, independently of $T$:
 \begin{equation}\label{p002}
  \int_{\R^d} |\nabla \sigma_{T,R}|^r \le C(d,r) \int_{\R^d} |\eta_R f|^r \le C(d,r) \int_{B_{2R}} |f|^r.
 \end{equation}
 
 We fix $T > 0$, and for $R' \ge R \gg \sqrt{T}$ we consider the difference $\sigma_{T,R}(x) - \sigma_{T,R'}(x)$ for points $x \in B_{R/2}(0)$. From the pointwise estimates on the massive Green's function $G_T$ (see, e.g.,~\cite[Corollary 1.5]{gloriamarahrens}) of the form
 \begin{align*}
  |\nabla G_T(x,y)| &\le C\frac{e^{-c\frac{1}{\sqrt{T}} |x-y|}}{|x-y|^{d-1}},
 \\
  |\nabla_x\nabla_y G_T(x,y)| &\le C\frac{e^{-c\frac{1}{\sqrt{T}} |x-y|}}{|x-y|^{d}},
 \end{align*}
 we get that for $x \in B_{R/2}(0)$ 
 \begin{equation}\label{p001}
  R^{-1} |(\sigma_{T,R}-\sigma_{T,R'})(x)| + |\nabla(\sigma_{T,R}-\sigma_{T,R'})(x)| \le C \frac{\sqrt{T}}{R} e^{-cR/\sqrt{T}}\|f\|_{L^\infty} \le C e^{-cR/\sqrt{T}\|f\|_{L^\infty}}. 
 \end{equation}
 In particular, in the limit $R \to \infty$ we have that $\sigma_{T,R}$ converges (pointwise) to $\sigma_T$, where $\sigma_T$ is stationary and satisfies 
 \begin{equation}\nonumber
  \tfrac{1}{T} \sigma_T - \Delta \sigma_T = -\nabla \cdot f.
 \end{equation}
 Moreover, estimate \eqref{p001} with $\sigma_{T,R'}$ replaced by $\sigma_T$ (estimate \eqref{p001} does not depend on $R'$, and so we are allowed to perform the limit $R' \to \infty$) in particular implies
 \begin{equation}\nonumber
  \int_{B_{R/2}(0)} |\nabla \sigma_{T,R} - \nabla \sigma_T|^r \le C R^d e^{-crR/\sqrt{T}} \|f\|^r_{L^\infty}. 
 \end{equation}
 We combine this estimate with \eqref{p002} to arrive at
 \begin{equation}\nonumber
  \fint_{B_{R/2}(0)} |\nabla \sigma_{T}|^r \le C e^{-crR/\sqrt{T}} \|f\|^r_{L^\infty} + C \fint_{B_{2R}(0)} |f|^r.
 \end{equation}
 Then by the ergodic theorem, as $R \to \infty$, the left-hand side converges to $\en{ |\nabla \sigma_T|^r }$ while the right-hand side converges to $\en{ |f|^r }$, i.e., we obtain
 \begin{equation}\nonumber
  \en{ |\nabla \sigma_T|^r } \le C \en{ |f|^r }. 
 \end{equation}
 Finally, since the sequence $\nabla \sigma_T$ is bounded in $L^r(\Omega)$, we can send $T \to \infty$ and obtain in the limit $\nabla \sigma$ which satisfies~\eqref{p003} and \eqref{p004}.

It remains to establish the ellipticity of the homogenized coefficient field $\ah$.  For the lower bound, we observe that, for an arbitrary $\xi\in\mathbb{R}^d$,
$$\xi\cdot \ah \xi=\en{\xi\cdot a(\nabla \phi_\xi+\xi)}=\en{(\nabla\phi_\xi+\xi)\cdot a(\nabla\phi_\xi+\xi)},$$
where the final inequality follows from the definition of the corrector $\phi_\xi$.  Therefore, by the definition of $\lambda$ in~\eqref{apq} we get
\begin{equation}\nonumber
\xi\cdot\ah\xi\geq\en{\lambda\abs{\nabla\phi_\xi+\xi}^2}\geq \en{\lambda^{-1}}^{-1}\abs{\xi}^2\geq\frac{1}{K}\abs{\xi}^2,
\end{equation}
where the last but one inequality follows from Jensen's inequality used for a jointly convex function $(f,g) \mapsto f^2/g$ with the choice $(f,g) = (\nabla \phi_\xi + \xi, \lambda^{-1})$, and the fact that $\en{\nabla\phi_\xi}=0$.
For the upper bound, for an arbitrary $\xi\in\mathbb{R}^d$, using the definition of $\mu$ from (\ref{apq}),
$$\abs{\ah \xi}=\abs{\en{a(\nabla\phi_\xi+\xi)}}\leq \en{\abs{a(\nabla\phi_\xi+\xi)}}\leq \en{\mu^{\frac{1}{2}}((\nabla\phi_\xi+\xi)\cdot a(\nabla\phi_\xi+\xi))^{\frac{1}{2}}}.$$
Then, after an application of H\"older's inequality, the definition of the corrector $\phi_\xi$ implies that
$$\abs{\ah \xi}\leq \en{\mu}^{\frac{1}{2}}\en{(\nabla\phi_\xi+\xi)\cdot a(\nabla\phi_\xi+\xi)}^{\frac{1}{2}}=\en{\mu}^{\frac{1}{2}}\en{\xi\cdot a(\nabla\phi_\xi+\xi)}^{\frac{1}{2}}\leq K^\frac{1}{2}\abs{\xi}^\frac{1}{2}\abs{\ah\xi}^{\frac{1}{2}}.$$
Dividing by $\abs{\ah\xi}^{\frac{1}{2}}$ yields the desired upper bound, and completes the proof.

\section{Proof of Lemma~\ref{lm2}}

To prove the sublinearity of the correctors $\phi$ and $\sigma$ we will only use that their gradients are stationary fields with zero expectation and that they have bounded $\frac{2q}{q+1}$ and $\frac{2p}{p+1}$ moments, respectively. Hence, we will only show the argument for $\phi$, the argument for $\sigma$ being analogous (after swapping $p$ and $q$). 

Concerning the corrector $\phi$, it is our aim to prove that \begin{equation}\label{phisublinear}  
\lim_{R\rightarrow\infty} 
\frac{1}{R} \biggl( \fint_{B_R} |\phi - \fint_{B_R} \phi|^{\frac{2p}{p-1}} \biggr)^{\frac{p-1}{2p}}=0.
\end{equation}  
Our proof is a simplified version of the proof a similar, seemingly slightly stronger, property (see ~\cite[Lemma 5.1]{ChiariniDeuschel2014}):
\begin{equation}\label{phisub2}
\lim_{R\rightarrow\infty} \frac{1}{R} \biggl( \fint_{B_R} |\phi|^{\frac{2p}{p-1}} \biggr)^{\frac{p-1}{2p}}=0.
\end{equation}
Before we prove~\eqref{phisublinear}, we point out that in fact it is equivalent~\eqref{phisub2}. Indeed, assuming~\eqref{phisublinear} for any $\delta > 0$ we find $r_0 > 0$ such that for all $R \ge r_0$
\begin{equation}\nonumber
\biggl( \fint_{B_R} |\phi - \fint_{B_R} \phi|^{\frac{2p}{p-1}} \biggr)^{\frac{p-1}{2p}} \le \delta R,
\end{equation}
from where by the triangle inequality we get for any $R \ge r_0$ and $R' \in [R,2R]$
\begin{equation}\nonumber
\biggl| \fint_{B_R} \phi - \fint_{B_{R'}} \phi \biggr| \le \biggl( \fint_{B_R} |\phi - \fint_{B_R} \phi|^{\frac{2p}{p-1}} \biggr)^{\frac{p-1}{2p}} + \biggl( \fint_{B_R} |\phi - \fint_{B_{R'}} \phi|^{\frac{2p}{p-1}} \biggr)^{\frac{p-1}{2p}} \le C \delta R.
\end{equation}
Hence by the dyadic argument we see that $\bigl| \fint_{B_{r_0}} \phi - \fint_{B_{R}} \phi \bigr| \le C \delta R$, which implies for $R \ge r_0$
\begin{equation}\nonumber
\frac{1}{R} \biggl( \fint_{B_R} |\phi|^{\frac{2p}{p-1}} \biggr)^{\frac{p-1}{2p}} \le
\frac{1}{R} \biggl( \fint_{B_R} |\phi - \fint_{B_R} \phi|^{\frac{2p}{p-1}} \biggr)^{\frac{p-1}{2p}} + C\delta + \frac{1}{R} \biggl| \fint_{B_{r_0}} \phi \biggr|,
\end{equation}
from where we get that $\limsup_{R\to \infty} \frac{1}{R} \left( \fint_{B_R} |\phi|^{\frac{2p}{p-1}} \right)^{\frac{p-1}{2p}} \le C\delta$, and~\eqref{phisub2} immediately follows. 

Let us now show the argument for~\eqref{phisublinear}, which is essentially an immediate consequence of the ergodic theorem, the Sobolev/Rellich-Kondrachov embedding and our assumption \begin{equation}\label{phif2} \frac{1}{p}+\frac{1}{q}  < \frac{2}{d}.\end{equation}  Fix $i\in\{1,\ldots,d\}$ and consider, for each $\epsilon\in(0,1)$, the rescaling $\phi^\epsilon_i(\cdot)=\epsilon\phi_i(\frac{\cdot}{\epsilon})$.  Assumption (\ref{phif2}) and the Sobolev embedding theorem imply that, for each $\epsilon\in(0,1)$, 
\begin{equation}\label{phisublinear1}
\biggl( \int_{B_1} |\phi^\epsilon_i - \fint_{B_1} \phi_i^\epsilon|^{\frac{2p}{p-1}} \biggr)^{\frac{p-1}{2p}}\lesssim \norm{\nabla \phi_i^\epsilon}_{L^{\frac{2q}{q+1}}(B_1)}.
\end{equation}  
Since the estimates contained in (\ref{c1}) and the ergodic theorem coupled with the stationarity and ergodicity of the environment imply that, for $\en{\cdot}$-a.e $a$, the gradient $\nabla \phi_i^\epsilon$ converges weakly to zero in $L^{\frac{2q}{q+1}}(B_1)$, we have for the renormalizations
$$(\phi_i^\epsilon-\fint_{B_1}\phi^\epsilon_i)\rightharpoonup 0 \;\; \textrm{weakly in}\;\; W^{1,\frac{2q}{q+1}}(B_1).$$
Finally, since the weak convergence and (\ref{phisublinear1}) imply that, for $\en{\cdot}$-a.e. environment, the sequence $\{(\phi^\epsilon_i-\fint_{B_1}\phi^\epsilon_i)\}_{\epsilon\in(0,1)}$ is bounded in $W^{1,\frac{2q}{q+1}}(B_1)$, the compactness of the embedding $W^{1,\frac{2q}{q+1}}(B_1)\hookrightarrow L^{\frac{2p}{p-1}}(B_1)$, owing to the strict inequality in (\ref{phif2}), implies for $\en{\cdot}$-a.e. $a$, the strong convergence
$$0 = \lim_{\epsilon\rightarrow 0}\biggl( \int_{B_1} |\phi^\epsilon_i - \fint_{B_1} \phi_i^\epsilon|^{\frac{2p}{p-1}} \biggr)^{\frac{p-1}{2p}}=\lim_{R\rightarrow 0}\frac{1}{R}\biggl( \fint_{B_R} |\phi_i - \fint_{B_R} \phi_i|^{\frac{2p}{p-1}} \biggr)^{\frac{p-1}{2p}}.$$
This, since $i\in\{1,\ldots,d\}$ was arbitrary, completes the argument for $\phi$.
 
\section{Proof of Theorem~\ref{prop1}}
 
 The strategy of the proof of the theorem is very similar to the proof of a similar proposition in~\cite{GNO4}. The idea is to first show decay of excess for one value $\theta_0$ of the ratio $r/R$, and then iterate this estimate to show excess decay for all values $r/R$.
 
 To show the decay for a fixed value of $r/R$, the idea is to estimate the homogenization error in $B_R$ determined by the difference between the $a$-harmonic function $u$ and a correction of an appropriately chosen $\ah$-harmonic function to be denoted $v$.  In the uniformly elliptic setting, following the arguments of~\cite{GNO4}, the boundary values of the $\ah$-harmonic function can be chosen to coincide with $u$ on a sphere with generic radius close to $R$.  In the non-uniformly elliptic case, it is necessary, as explained in Step 2 of our arguments, to consider a $v$ which agrees on a generic sphere with $u_\epsilon$ in the ``Dirichlet Case'' $q\geq p$ and which satisfies $\nu \cdot \ah \nabla v = (\nu\cdot a\nabla u)_\epsilon$ in the ``Neumann Case'' $p\geq q$, where the subscript $\epsilon$ denotes the convolution at scale $\epsilon$ with a smoothing kernel on the sphere.
 
The corresponding augmented homogenization error will be defined for an appropriately chosen cutoff function $\eta$ in $B_R$ as $w:=u-(1+\eta\phi_i\partial_i )v$. In Step 1, as in \cite{GNO4}, we derive the equation satisfied by $w$, and in Step 2 use this equation to obtain energy estimates for the homogenization error, without the cutoff $\eta$, on a smaller ball.  The argument is concluded in Steps 3, 4 and 5, where the iterative argument of \cite{GNO4} is used to obtain the statement on excess decay.

 \medskip
  {\bf Step 1.}  Let $u$ be an $a$-harmonic function in $B_1$.  In this step we consider the augmented homogenization error
\begin{equation}\nonumber
w:=u-(1+\eta\phi_i\partial_i)v,
\end{equation}
defined by a smooth function $\eta$ and an $\ah$-harmonic function $v$ in $B_1$.  It will be shown now that $w$ solves the divergence-form equation
\begin{equation}\label{errorequation}
-\nabla\cdot a\nabla w = \nabla\cdot \left((1-\eta)(a-\ah)\nabla v+(\phi_ia-\sigma_i)\nabla(\eta\partial_iv)\right)\;\;\textrm{in}\;\;B_1,
\end{equation}
where the crucial ingredient of the proof is the skew-symmetric flux correction $\sigma$.  We remark that the flux correction was used previously in the context of periodic homogenization (see, e.g.,~\cite{ZikovKozlovOlejnik}), and in stochastic homogenization it was introduced only recently in~\cite{GNO4}, where ~\eqref{errorequation} was first derived.

For convenience of the reader we repeat here the computation leading to~\eqref{errorequation}, and to keep the notation lean, in this and the following steps, we will without loss of generality assume that the components of $\phi$ and $\sigma$ have zero spatial average on $B_1$.  Otherwise, for each $i\in\{1,\ldots,d\}$, we would replace $\phi_i$ with $(\phi_i-\fint_{B_1}\phi_i)$, and similarly for $\sigma$.

First, compute the gradient of $w$ to find
\begin{equation}\nonumber
\nabla w=\nabla u-(\nabla v+\eta\partial_iv\nabla \phi_i+\phi_i\nabla(\eta\partial_iv)),
\end{equation}
then we use the $a$-harmonicity of $u$ to obtain
\begin{equation}\nonumber
-\nabla\cdot a\nabla w=\nabla\cdot a\nabla v+\nabla\cdot a(\eta\partial_iv\nabla \phi_i+\phi_i\nabla(\eta\partial_iv)).
\end{equation}
Since
$$\nabla\cdot a(\eta\partial_iv\nabla \phi_i)=\nabla\cdot (\eta\partial_iv a(\nabla \phi_i+e_i))-\nabla \cdot \eta a\nabla v,$$
the vanishing divergence $-\nabla\cdot a(\nabla \phi_i+e_i)=0$ implies that
\begin{equation}\nonumber
-\nabla\cdot a\nabla w=\nabla\cdot (1-\eta)a\nabla v+\nabla(\eta\partial_iv)\cdot a(\nabla\phi_i+e_i)+\nabla\cdot(\phi_ia\nabla(\eta\partial_iv)).
\end{equation}
Then, after observing both that
\begin{equation}\nonumber
\ah e_i\cdot\nabla(\eta\partial_iv)=\nabla\cdot(\eta\partial_iv \ah e_i)=\nabla\cdot(\eta \ah\nabla v),
\end{equation}
and, since $-\nabla \cdot\ah\nabla v=0$, that
$$\nabla\cdot (\eta \ah\nabla v)=-\nabla\cdot (1-\eta)\nabla v,$$ we have $$-\nabla\cdot a\nabla w=\nabla\cdot((1-\eta)(a-\ah )\nabla v)+\nabla(\eta\partial_iv)\cdot(a(\nabla\phi_i+e_i)-\ah e_i)+\nabla\cdot(\phi_ia\nabla(\eta\partial_iv)).$$

The skew-symmetry of the flux correction $\sigma$ now plays a role.  Since $$\nabla\cdot \sigma_i=q_i=a(\phi_i+e_i)-\ah e_i,$$ we have, for an arbitrary test function $\zeta$, the distributional identity $$\nabla \zeta\cdot (\nabla\cdot\sigma_i)=\partial_j\zeta\partial_k\sigma_{ijk}=\partial_k(\partial_j\zeta\sigma_{ijk})=\partial_k(\sigma_{ijk}\partial_j\zeta)=-\nabla\cdot(\sigma_i\nabla \zeta),$$ from which~\eqref{errorequation} follows. This completes the proof of this step.

\medskip

{\bf{Step 2.}}  The boundary conditions for $v$ and the cutoff $\eta$ are now specified in order to use equation \eqref{errorequation} to obtain an energy estimate for the homogenization error.  We remark that the arguments will be carried out for the unit ball $B_1$, and the general statement will be obtained by scaling.  We assume that 
\begin{equation}\label{f1}
\big(\int_{B_1}\lambda^{-q}\big)^\frac{1}{q}+\big(\int_{B_1}\mu^p\big)^\frac{1}{p} \le \Lambda, \qquad \textrm{where } \frac{1}{p} + \frac{1}{q} \le \frac{2}{d},
\end{equation}
with
\begin{equation}\label{f2}
\lambda:=\inf_{\xi\in\mathbb{R}^d}\frac{a \xi\cdot \xi}{|\xi|^2}\quad\mbox{and}\quad
\mu:=\sup_{\xi\in\mathbb{R}^d}\frac{|a\xi|^2}{a \xi\cdot \xi},
\end{equation}
and consider an $a$-harmonic function $u$ in $B_1$, that is,
\begin{equation}\label{f7}
-\nabla\cdot a\nabla u=0\quad\mbox{in}\;B_1,
\end{equation}
where by homogeneity we may assume that
\begin{equation}\label{f4}
\int_{B_1}\nabla u\cdot a\nabla u=1.
\end{equation}
We will construct an $\ah$-harmonic function $v$ in $B_\frac{1}{2}$ satisfying, in view of the normalization (\ref{f4}),
\begin{equation}\label{36+}
\int_{B_\frac{1}{2}}\nabla v\cdot \ah \nabla v\lesssim \Lambda \int_{B_1}\nabla u\cdot a\nabla u=\Lambda,
\end{equation}
and for which the homogenization error $w:=u-(1+\phi_i\partial_i)v$ satisfies
\begin{align}\label{f50}
\int_{B_\frac{1}{4}}\nabla w\cdot a \nabla w &\lesssim \Lambda\epsilon^{1-(\frac{1}{2p}+\frac{1}{2q})(d-1)} +\Lambda^2\rho^{\min\{\frac{p-1}{2p},\frac{q-1}{2q}\}}\frac{1}{\epsilon^{d\min\{\frac{p+1}{p},\frac{q+1}{q}\}}}\nonumber\\
&\quad +\Lambda^2\frac{1}{\rho^{d+2}}
\Big(\big(\int_{B_1}|\phi|^{\frac{2p}{p-1}}\big)^\frac{p-1}{p}+
\big(\int_{B_1}|\sigma|^{\frac{2q}{q-1}}\big)^\frac{q-1}{q}\Big),
\end{align}
for any fudge factors $\epsilon \in (0,1]$ and $\rho \in (0,\frac{1}{8})$.
We recall that $\lesssim$ denotes a constant depending only upon the dimension $d$ and the constant $K$ from (\ref{apq}) through the ellipticity of the homogenized coefficients.

We begin now with the construction of $v$ which will in fact be an $\ah$-harmonic on a somewhat larger ball $B_r$, for some suitably chosen radius $r\in[\frac{1}{2},1]$:
\begin{equation}\label{f6}
-\nabla\cdot \ah\nabla v=0\quad\mbox{on}\;B_r.
\end{equation}
The idea is to distinguish the two cases
\begin{equation}\label{f5}
\mbox{the ``Dirichlet case''}\;q\ge p\quad\mbox{and the}\quad\mbox{``Neumann case''}\;p\ge q.
\end{equation}
In the Dirichlet case, we define $v$ via the Dirichlet boundary condition
\begin{equation}\label{f23}
v=u_\epsilon\;\mbox{on}\;\partial B_r,
\end{equation}
whereas in the Neumann case, we impose
\begin{equation}\label{f19}
\nu\cdot \ah\nabla v=(\nu\cdot a\nabla u)_\epsilon\;\mbox{on}\;\partial B_r.
\end{equation}
Here the subscript $\epsilon$ stands for a convolution on $\partial B_r$ with scale $\epsilon>0$.  

Since H\"older's inequality, (\ref{f1}), (\ref{f2}), and (\ref{f4}) imply
\begin{align*}
\biggl(\int_{B_1}|\nabla u|^\frac{2q}{q+1}\biggr)^\frac{q+1}{2q} + \biggl(\int_{B_1}|a\nabla u|^\frac{2p}{p+1}\biggr)^\frac{p+1}{2p} 
\le\biggl( \biggl(\int_{B_1}\lambda^{-q}\biggr)^{\frac{1}{2q}} + \biggl(\int_{B_1}\mu^{p}\biggr)^{\frac{1}{2p}} \biggr)
\biggl(\int_{B_1}\nabla u\cdot a\nabla u\biggr)^\frac{1}{2}
\lesssim \Lambda^\frac{1}{2},
\end{align*}
we can find a radius $r \in [\frac{1}{2},1]$ such that both the field and the current of $u$ have
the same integrability on $\partial B_r$ as on $B_1$, in the sense that
\begin{equation}\label{f13}
  \biggl(\int_{\partial B_r}|\nabla u|^\frac{2q}{q+1}\biggr)^\frac{q+1}{2q} + \biggl(\int_{\partial B_r}|a\nabla u|^\frac{2p}{p+1}\biggr)^\frac{p+1}{2p} \lesssim \Lambda^{\frac{1}{2}}. 
\end{equation}
Using that both estimates are preserved by convolution, it follows that
\begin{equation}\nonumber
\biggl(\int_{\partial B_r}|\nabla^{\textrm{tan}} v|^\frac{2q}{q+1}\biggr)^\frac{q+1}{2q}
\lesssim\Lambda^\frac{1}{2}\quad\mbox{in the Dirichlet case}
\end{equation}
and
\begin{equation}\nonumber
\biggl(\int_{\partial B_r}|\nu\cdot \ah\nabla v|^\frac{2p}{p+1}\biggr)^\frac{p+1}{2p}
\lesssim
\Lambda^\frac{1}{2}\quad\mbox{in the Neumann case}.
\end{equation}
By constant-coefficient elliptic theory applied to the Dirichlet-to-Neumann operator (see, e.g. Stein \cite[Chapter~7]{SteinHarmonic}) this yields
\begin{equation}\label{f10+}
\big(\int_{\partial B_r}|\nabla v|^\frac{2q}{q+1}\big)^\frac{q+1}{2q}
\lesssim
\Lambda^\frac{1}{2}\quad\mbox{in the Dirichlet case}
\end{equation}
and 
\begin{equation}\label{f20}
\big(\int_{\partial B_r}|\nabla v|^\frac{2p}{p+1}\big)^\frac{p+1}{2p}
\lesssim
\Lambda^\frac{1}{2}\quad\mbox{in the Neumann case}.
\end{equation}
These estimates motivate the case distinction (\ref{f5}), which we now use to prove (\ref{36+}).  
To simplify the notation, we assume without loss of generality that $v$ has zero average on $\partial B_r$, and test (\ref{f6}) with $v$ to obtain
$$\int_{B_r}\nabla v\cdot\ah\nabla v=-\int_{\partial B_r}v(\nu\cdot\ah\nabla v).$$
Using H\"older's inequality, the Sobolev embedding theorem and assumption (\ref{apq}), we have
$$\int_{B_r}\nabla v\cdot \ah\nabla v\lesssim \big(\int_{\partial B_r}|\nabla v|^\frac{2q}{q+1}\big)^\frac{q+1}{2q}\big(\int_{\partial B_r}|v|^\frac{2q}{q-1}\big)^\frac{q-1}{2q}\lesssim \big(\int_{\partial B_r}|\nabla v|^\frac{2q}{q+1}\big)^\frac{q+1}{2q}\big(\int_{\partial B_r}|\nabla v|^\frac{2p}{p+1}\big)^\frac{p+1}{2p}.$$
Therefore, in view of (\ref{f10+}), (\ref{f20}) and the case distinction (\ref{f5}), we conclude that
\begin{equation}\label{f36++}\int_{B_r}\nabla v\cdot\ah\nabla v\lesssim\Lambda,\end{equation}
which implies (\ref{36+}) since $r\in[\frac{1}{2},1]$.

We now specify precisely the cutoff function defining the augmented homogenization error.  Let $0<\rho<\frac{1}{8}$ to be fixed later and let $0\leq \eta \leq 1$ be a smooth function satisfying
\begin{equation}\label{ee_cutoff} 
 \eta=\left\{\begin{array}{lll} 1 & \textrm{on}\;\;\overline{B}_{r-2\rho}, \\ 0 & \textrm{on}\;\;B_r\setminus B_{r-\rho},\end{array}\right.\; \textrm{with}\;\;\abs{\nabla\eta}\lesssim \frac{1}{\rho}.
\end{equation}
We recall from Step 1, see (\ref{errorequation}), the formula for the augmented homogenization error
\begin{equation}\nonumber
-\nabla\cdot a\nabla w=-\nabla\cdot\big((1-\eta)(\ah-a)\nabla v\big)
+\nabla\cdot\big((\phi_ia-\sigma_i)\nabla(\eta\partial_iv)\big)\;\;\textrm{in}\;\;B_{r},
\end{equation}
which requires both (\ref{f7}) and (\ref{f6}), and test this equation with $w$ to obtain
\begin{align}
\int_{B_r}\nabla w\cdot a\nabla w
&=\int_{\partial B_r}(u-v)\nu\cdot(a\nabla u-\ah\nabla v)\nonumber\\
&\quad +\int_{B_r}(1-\eta)\nabla w\cdot(\ah-a)\nabla v\nonumber\\
&\quad -\int_{B_r}\nabla w\cdot (\phi_ia-\sigma_i)\nabla(\eta\partial_iv).\nonumber
\end{align}
Appealing to (\ref{f2}) we obtain by Young's inequality
\begin{align*}
\int_{B_r}\nabla w\cdot a\nabla w
&\lesssim\big|\int_{\partial B_r}(u-v)\nu\cdot(a\nabla u-\ah\nabla v)\big|\\
&\quad+\int_{B_r}(1-\eta)^2(\mu+\frac{1}{\lambda})|\nabla v|^2\nonumber\\
&\quad+\int_{B_r}(\mu\phi_i^2+\frac{1}{\lambda}|\sigma_i|^2)|\nabla(\eta\partial_iv)|^2.\nonumber
\end{align*}
Then, using the properties of the cutoff $\eta$ from (\ref{ee_cutoff}),
\begin{align}
\int_{B_r}\nabla w\cdot a\nabla w
&\lesssim\big|\int_{\partial B_r}(u-v)\nu\cdot(a\nabla u-\ah\nabla v)\big|\nonumber\\
&\quad+\int_{B_r\setminus B_{r-2\rho}}(\mu+\frac{1}{\lambda})|\nabla v|^2\nonumber\\
&\quad+\int_{B_1}(\mu|\phi|^2+\frac{1}{\lambda}|\sigma|^2)\sup_{B_{r-\rho}}(|\nabla^2v|+\frac{1}{\rho}|\nabla v|)^2.\nonumber
\end{align}
In view of (\ref{f1}), this yields by H\"older's inequality
\begin{align}
\lefteqn{\int_{B_r}\nabla w\cdot a\nabla w}\nonumber\\
&\lesssim\big|\int_{\partial B_r}(u-v)\nu\cdot(a\nabla u-\ah\nabla v)\big|\label{f18}\\
&\quad+\Lambda\rho^{\min\{\frac{p-1}{2p},\frac{q-1}{2q}\}}(\int_{B_r}|\nabla v|^{\max\{\frac{4p}{p-1},\frac{4q}{q-1}\}})^{\min\{\frac{p-1}{2p},\frac{q-1}{2q}\}}\label{f12}\\
&\quad+\Lambda\Big(\big(\int_{B_1}|\phi|^{\frac{2p}{p-1}}\big)^\frac{p-1}{p}+
\big(\int_{B_1}|\sigma|^{\frac{2q}{q-1}}\big)^\frac{q-1}{q}\Big)\sup_{B_{r-\rho}}(|\nabla^2v|+\frac{1}{\rho}|\nabla v|)^2.\label{f17}
\end{align}

For the last term, we use multiple times the Caccioppoli inequality for uniformly elliptic systems, see \cite[Theorem~7.1]{giaqmart}, the uniform ellipticity of the homogenized coefficient field, the Sobolev embedding theorem, and (\ref{f6}) to obtain
\begin{equation}\nonumber
\sup_{B_{r-\rho}}(|\nabla^2v|+\frac{1}{\rho}|\nabla v|)^2
\lesssim \frac{1}{\rho^{d+2}}\int_{B_r}\nabla v\cdot \ah \nabla v,
\end{equation}
where the full argument can be found in the proof of \cite[Lemma~3]{GNO4}.  Therefore, in view of (\ref{f36++}), for both the Dirichlet and Neumann cases we have
\begin{equation}\label{f16}
\sup_{B_{r-\rho}}(|\nabla^2v|+\frac{1}{\rho}|\nabla v|)^2
\lesssim \frac{1}{\rho^{d+2}}\Lambda.
\end{equation}

We now turn to the middle right-hand side term \eqref{f12}.  It follows from the $L^p$-theory for constant-coefficient elliptic systems, see \cite[Theorem~7.1]{giaqmart}, and an explicit radial extension of the smooth boundary data $u_\epsilon$ into the ball $B_r$ that, in the Dirichlet case, 
\begin{equation}\nonumber (\int_{B_r}|\nabla v|^{\max\{\frac{4p}{p-1},\frac{4q}{q-1}\}})^{\min\{\frac{p-1}{2p},\frac{q-1}{2q}\}}=(\int_{B_r}|\nabla v|^{\frac{4q}{q-1}})^{\frac{q-1}{2q}}\lesssim (\int_{\partial B_r} |\nabla^{\textrm{tan}} u_\epsilon|^{\frac{4q}{q-1}})^{\frac{q-1}{2q}}\lesssim (\sup_{\partial B_r}|\nabla^{\textrm{tan}} u_\epsilon|)^2.\end{equation}
We then use the inverse estimate
\begin{equation}\nonumber
\sup_{\partial B_r}|\nabla^{\textrm{tan}} u_\epsilon|\lesssim\frac{1}{\epsilon^{d\frac{q+1}{2q}}}
\big(\int_{\partial B_r}|\nabla u|^\frac{2q}{q+1}\big)^\frac{q+1}{2q}\quad\mbox{in the Dirichlet case},
\end{equation}
so that together with (\ref{f13}) we obtain
\begin{equation}\nonumber
(\int_{B_r}|\nabla v|^{\max\{\frac{4p}{p-1},\frac{4q}{q-1}\}})^{\min\{\frac{p-1}{2p},\frac{q-1}{2q}\}}\lesssim\frac{1}{\epsilon^{d\frac{q+1}{q}}}\Lambda\quad\mbox{in the Dirichlet case}.
\end{equation}
The analogous estimates for the the Neumann case yield, with help of~\eqref{f13}, 
\begin{equation}\nonumber
(\int_{B_r}|\nabla v|^{\max\{\frac{4p}{p-1},\frac{4q}{q-1}\}})^{\min\{\frac{p-1}{2p},\frac{q-1}{2q}\}}=(\int_{B_r}|\nabla v|^{\frac{4p}{p-1}})^{\frac{p-1}{2p}}\lesssim\frac{1}{\epsilon^{d\frac{p+1}{p}}}\Lambda,
\end{equation}
so that, in combination, both cases satisfy
\begin{equation}\label{f15}
(\int_{B_r}|\nabla v|^{\max\{\frac{4p}{p-1},\frac{4q}{q-1}\}})^{\min\{\frac{p-1}{2p},\frac{q-1}{2q}\}}\lesssim\frac{1}{\epsilon^{d\min\{\frac{p+1}{p},\frac{q+1}{q}\}}}
\Lambda.
\end{equation}

It remains to treat the boundary term (\ref{f18}). We 
first treat the Neumann case (\ref{f19}), for which
we may appeal to the symmetry of the convolution operator to write
\begin{equation}\nonumber
\int_{\partial B_r}(u-v)\nu\cdot(a\nabla u-\ah\nabla v)
=\int_{\partial B_r}((u-v)-(u-v)_\epsilon)\nu\cdot a\nabla u,
\end{equation}
so that we obtain by H\"older's inequality together with~\eqref{f13}
\begin{align}\label{f21bis}
\lefteqn{\big|\int_{\partial B_r}(u-v)\nu\cdot(a\nabla u-\ah\nabla v)\big|}\nonumber\\
&\lesssim\Lambda^\frac{1}{2}\Big(\big(\int_{\partial B_r}|u-u_\epsilon|^\frac{2p}{p-1}\big)^\frac{p-1}{2p}+
\big(\int_{\partial B_r}|v-v_\epsilon|^\frac{2p}{p-1}\big)^\frac{p-1}{2p}\Big).
\end{align}
Interpolating between a convolution estimate of the form
\begin{equation}\nonumber
\big(\int_{\partial B_r}|u-u_\epsilon|^\frac{2p}{p-1}\big)^\frac{p-1}{2p}
\lesssim\epsilon \big(\int_{\partial B_r}|\nabla^{\textrm{tan}} u|^\frac{2p}{p-1}\big)^\frac{p-1}{2p}
\end{equation}
and an implication of the triangle and the Sobolev inequalities
\begin{equation}\nonumber
\big(\int_{\partial B_r}|u-u_\epsilon|^\frac{2p}{p-1}\big)^\frac{p-1}{2p}
\le 
\big(\int_{\partial B_r}|u|^\frac{2p}{p-1}\big)^\frac{p-1}{2p} + 
\big(\int_{\partial B_r}|u_\epsilon|^\frac{2p}{p-1}\big)^\frac{p-1}{2p} 
\lesssim 
\big(\int_{\partial B_r}|\nabla^{\textrm{tan}} u|^s\big)^\frac{1}{s},
\end{equation}
where $\frac{1}{s} = \frac{p-1}{2p} + \frac{1}{d-1}$, yields
\begin{equation}\nonumber
 \big(\int_{\partial B_r}|u-u_\epsilon|^\frac{2p}{p-1}\big)^\frac{p-1}{2p}
\lesssim\epsilon^{1-(\frac{1}{2p}+\frac{1}{2q})(d-1)}
\big(\int_{\partial B_r}|\nabla^{\textrm{tan}} u|^\frac{2q}{q+1}\big)^\frac{q+1}{2q}.
\end{equation}
The analogous estimate for $v$ reads 
\begin{equation}
\big(\int_{\partial B_r}|v-v_\epsilon|^\frac{2p}{p-1}\big)^\frac{p-1}{2p}
\lesssim\epsilon^{1-\frac{1}{p}(d-1)}
\big(\int_{\partial B_r}|\nabla^{\textrm{tan}} v|^\frac{2p}{p+1}\big)^\frac{p+1}{2p}.\nonumber
\end{equation}
We plug~\eqref{f13} and~\eqref{f20} into these two estimates and use our case distinction (\ref{f5}) to arrive at
\begin{align}\nonumber
\big(\int_{\partial B_r}|u-u_\epsilon|^\frac{2p}{p-1}\big)^\frac{p-1}{2p}
+\big(\int_{\partial B_r}|v-v_\epsilon|^\frac{2p}{p-1}\big)^\frac{p-1}{2p}
\lesssim\Lambda^\frac{1}{2}\epsilon^{1-(\frac{1}{2p}+\frac{1}{2q})(d-1)}.
\end{align}
Inserting this into (\ref{f21bis}) we obtain for the boundary term in the Neumann case
\begin{align}\label{f25}
\big|\int_{\partial B_r}(u-v)\nu\cdot(a\nabla u-\ah\nabla v)\big|
\lesssim\Lambda\epsilon^{1-(\frac{1}{2p}+\frac{1}{2q})(d-1)}.
\end{align}

\medskip

We finally turn to the Dirichlet case (\ref{f23}) for the boundary term (\ref{f18}), which yields a simpler estimate as compared to the Neumann case because $u$ and $v$ are immediately comparable along the boundary $\partial B_r$.  We use H\"older's inequality, the triangle inequality, and the case distinction in combination with~\eqref{f13} 
and (\ref{f10+}) to obtain
\begin{align}\nonumber
\big|\int_{\partial B_r}(u-v)\nu\cdot(a\nabla u-\ah\nabla v)\big|
&\lesssim\Lambda^\frac{1}{2}\big(\int_{\partial B_r}|u-u_\epsilon|^\frac{2q}{q-1}\big)^\frac{q-1}{2q}.
\end{align}
Appealing again to convolution estimates used to obtain (\ref{f25}) and the case distinction, we obtain with (\ref{f13}) the estimate
\begin{equation}\label{f25+}
\big|\int_{\partial B_r}(u-v)\nu\cdot(a\nabla u-\ah\nabla v)\big| \lesssim \Lambda^{\frac{1}{2}}\epsilon^{1-\frac{1}{q}(d-1)}\big(\int_{\partial B_r}|\nabla^{\textrm{tan}} u|^\frac{2q}{q+1}\big)^\frac{q+1}{2q}\lesssim \Lambda \epsilon^{1-(\frac{1}{2p}+\frac{1}{2q})(d-1)},\end{equation}
which corresponds with the Neumann estimate (\ref{f25}). Inserting (\ref{f25+}) and (\ref{f25}) together with (\ref{f15}) and (\ref{f16})
into (\ref{f17}) yields (\ref{f50}).  This completes the proof of Step 2.

\medskip
{\bf Step 3.} In this step, we prove that whenever $u$ is an $a$-harmonic function on $B_R$, there exists for any $\delta_0>0$ a constant $C_0=C_0(\delta_0,d,\Lambda,K)>0$ such that, whenever on the ball $B_R$  the augmented corrector $(\phi,\sigma)$ satisfies (\ref{corrsmall}) with constant $C_0$, we have
\begin{equation}\label{f511}\fint_{B_{\frac{R}{4}}}\nabla w\cdot a \nabla w\leq \delta_0\fint_{B_R} \nabla u\cdot a\nabla u,\end{equation}
where $w:=u-(1+\phi_i\partial_i) v$ is the homogenization error.

Here comes the argument. 
Given $\delta_0>0$, owing to the linearity and scaling of the equation the estimate~\eqref{f50} from Step 2 implies existence of $C_1=C_1(d,\Lambda,K)$ such that
%
\begin{equation}\nonumber
\fint_{B_\frac{R}{4}}\nabla w\cdot a \nabla w \le C_1\Big(\Lambda\epsilon^{1-(\frac{1}{2p}+\frac{1}{2q})(d-1)} +\Lambda^2\rho^{\min\{\frac{p-1}{2p},\frac{q-1}{2q}\}}\frac{1}{\epsilon^{d\min\{\frac{p+1}{p},\frac{q+1}{q}\}}}+\Lambda^2\frac{1}{\rho^{d+2}C_0}\Big)\fint_{B_R}\nabla u\cdot a\nabla u.
\end{equation}
We first fix $\epsilon_0=\epsilon_0(\delta_0,d,\Lambda,K)>0$ small enough such that
$$C_1\Lambda\epsilon_0^{1-(\frac{1}{2p}+\frac{1}{2q})(d-1)}\leq \frac{1}{3}\delta_0,$$
which is possible in view of (\ref{pq}). Second, we choose $\rho_0=\rho_0(\delta_0,d,\Lambda,K)>0$ small enough satisfying
$$C_1\Lambda^2\rho_0^{\min\{\frac{p-1}{2p},\frac{q-1}{2q}\}}\frac{1}{\epsilon_0^{d\min\{\frac{p+1}{p},\frac{q+1}{q}\}}}\leq \frac{1}{3}\delta_0.$$
Finally, we select $C_0=C_0(\delta_0,d,\Lambda,K)>0$ large enough so that
$$C_1\Lambda^2\frac{1}{\rho_0^{d+2}C_0}\leq \frac{1}{3}\delta_0.$$
Since the right-hand sides in the three previous relations add to $\delta_0$, the proof of this step is complete. 

\medskip
{\bf Step 4.} 
For any $\alpha \in (0,1)$ there exists $\theta_0 \in (0,\frac{1}{4})$ and $C_0$ such that the following holds: For any radius $R > 0$ with the property that the augmented corrector is small on scale $R$:
 \begin{equation}
  \begin{aligned}
  \frac{1}{R} \biggl( \fint_{B_R} |\phi - \fint_{B_R} \phi|^{\frac{2p}{p-1}} \biggr)^{\frac{p-1}{2p}} &\le \frac{1}{C_0},
  \\
  \frac{1}{R} \biggl( \fint_{B_R} |\sigma - \fint_{B_R} \sigma|^{\frac{2q}{q-1}} \biggr)^{\frac{q-1}{2q}} &\le \frac{1}{C_0},
  \end{aligned}
 \end{equation}
and the corrector $\phi$ is small on scale $r := \theta_0 R$
 \begin{equation}
  \frac{1}{r} \biggl( \fint_{B_r} |\phi - \fint_{B_r} \phi|^{\frac{2p}{p-1}} \biggr)^{\frac{p-1}{2p}} \le \frac{1}{C_0},
 \end{equation}
we get that for every $a$-harmonic function $u$ on $B_R$ its excess satisfies
%
\begin{equation}\label{excessstep}
 \mathrm{Exc}(r) \le \theta_0^{2\alpha} \mathrm{Exc}(R),
\end{equation}
where we recall
$$\mathrm{Exc}(r)=\inf_{\xi\in\mathbb{R}^d}\fint_{B_r}(\nabla u-(\xi+\phi_\xi))\cdot a(\nabla u-(\xi+\phi_\xi)).$$

Here comes the argument. 
To simplify the notation, we will use $\| \cdot \|_{a,r}$ to denote the $L^2(a)$-intrinsic energy of vector fields~$U$ as defined by
\begin{equation}\nonumber
 \| U \|_{a,r} := \int_{B_r} U\cdot a U. 
\end{equation}
{\renewcommand{\norm}[2]{\| #1 \|_{a,#2}}%
We now use the definition of the homogenization error $w=u-(1+\phi_i\partial_i)v$, and observe that, for every $0<r\leq\frac{R}{4}$, we have
\begin{equation}\nonumber
 \norm{\nabla w}{r} = \norm{ \nabla u - \nabla v (\textrm{Id} + \nabla \phi) - \phi_i\nabla(\partial_i v)}{r}.
\end{equation}
The first two terms are decomposed as
$$\nabla u - \nabla v(\textrm{Id} + \nabla \phi) = \nabla u - \nabla v(0) (\textrm{Id} + \nabla \phi) + (\nabla v(0) - \nabla v) (\textrm{Id} + \nabla \phi)$$
 and, by defining $\xi := \nabla v(0)$, we obtain using the triangle inequality, since $0<r\leq\frac{R}{4}$,
\begin{equation}\label{f70}
\begin{aligned}
 \norm{ \nabla u - (\xi + \nabla \phi_\xi) }{r} &\leq \norm{\nabla w}{\frac{R}{4}} + (\sup_{B_r} |\nabla v(0) - \nabla v|)^2 (\norm{ \textrm{Id} + \nabla \phi }{r}) 
 \\
 &\quad+ (\sup_{B_r} |\nabla (\partial_iv)|)^2 \norm{\phi}{r}.
\end{aligned}
\end{equation}

The first term on the right-hand side of (\ref{f70}) is controlled using the estimate (\ref{f511}) from Step 3.  There exists, for any $\delta_0>0$, a constant $C_0=C_0(\delta_0,d,\Lambda,K)>0$ such that, whenever the corrector and flux correction satisfy (\ref{corrsmall}) with constant $C_0$, we have
\begin{equation}\label{f71}
\norm{\nabla w}{\frac{R}{4}}\le \delta_0\int_{B_R}\nabla u\cdot a\nabla u.
\end{equation}
It remains to control the last two terms on the right-hand side of (\ref{f70}).  First, since the corrector satisfies $$-\nabla \cdot a (\textrm{Id} + \nabla\phi) = 0,$$ the Caccioppoli estimate (Lemma~\ref{lmcac}) together with~\eqref{corrsmall} implies
 \begin{equation}\nonumber
  \norm{ \textrm{Id} + \nabla \phi }{r} \lesssim r^d.
 \end{equation}
Then, by repeating the argument leading to (\ref{f36++}), since $0<r\leq\frac{R}{4}$, we have the estimate
\begin{equation}\label{f74}
 \sup_{B_r} |\nabla^2 v|^2 \lesssim \frac{1}{R^{d+2}} \norm{ \nabla u }{R}.\end{equation}

We insert \eqref{f74} and \eqref{f71} into~\eqref{f70}, use (\ref{abound}) and (\ref{corrsmall}), and find for a constant $C_1=C_1(d,\Lambda, K)>0$ that
\begin{equation}\nonumber
 \norm{ \nabla u - (\xi + \nabla \phi_\xi) }{r} \leq C_1\biggl( \delta_0+ \left( \frac{r}{R} \right)^2 \left( \frac{r}{R} \right)^d\biggr)\norm{\nabla u}{R}.
\end{equation}
Then, after dividing by $r^d$,
\begin{equation}\label{f75}
 \fint_{B_r} (\nabla u - (\xi + \nabla \phi_\xi)) \cdot a(\nabla u - (\xi + \nabla \phi_\xi)) \le C_1 \biggl(\delta_0\Big(\frac{R}{r}\Big)^d+\Big(\frac{r}{R}\Big)^2\biggr)\fint_{B_R}  \nabla u \cdot a\nabla u.
 \end{equation}

Fix $\alpha\in(0,1)$. The proof of Step 4 will be complete once we prove that the ratio $\frac{r}{R}$ and the constant $\delta_0$ can fixed sufficiently small so as to guarantee the inequality
\begin{equation}\nonumber
C_1\biggl(\delta_0\Big(\frac{R}{r}\Big)^d+\Big(\frac{r}{R}\Big)^2\biggr) \le \left( \frac{r}{R} \right)^{2\alpha},
\end{equation}
and which is possible because $\alpha$ is strictly smaller than one.
First, we choose $\theta_0\in(0,\frac{1}{4})$ small enough to satisfy
\begin{equation}\label{f78}C_1\theta_0^2\leq \frac{1}{2}\theta_0^{2\alpha},\end{equation}
and then choose $\delta_0$ small enough to ensure
\begin{equation}\label{f79} C_1\delta_0\theta_0^d\leq \frac{1}{2}\theta_0^{2\alpha}.\end{equation}
Then, for $\theta_0$ satisfying (\ref{f78}) and (\ref{f79}), whenever $\frac{r}{R}=\theta_0$ we obtain from (\ref{f75}) the estimate
\begin{equation}\nonumber
 \inf_{\zeta \in \R^d} \fint_{B_r} (\nabla u - (\zeta + \nabla \phi_\zeta))\cdot a (\nabla u - (\zeta + \nabla \phi_\zeta)) \leq \left( \frac{r}{R} \right)^{2\alpha} \fint_{B_R} \nabla u \cdot a \nabla u.
 \end{equation}
After replacing $u$ with $u - \xi\cdot(x + \phi)$, where $\xi := \textrm{argmin}_{\xi \in \R^d} \fint_{B_R} a (\nabla u - (\xi + \nabla \phi_\xi)) \cdot (\nabla u - (\xi + \nabla \phi_\xi))$, we obtain~\eqref{excessstep} and complete Step 3.
}

\medskip
{\bf Step 5.} In the final step of the proof, we will prove that, for arbitrary pairs $r<R$ satisfying the hypothesis of Theorem~\ref{prop1}, we have the excess decay
\begin{equation}\label{excessstep1}
 \mathrm{Exc}(r) \lesssim \left( \frac{r}{R} \right)^{2\alpha} \mathrm{Exc}(R).
\end{equation}
For this, we iterate the estimate~\eqref{excessstep} to obtain, for $\theta_0\in(0,1)$ from Step 4,
\begin{equation}\nonumber
 \mathrm{Exc}(\theta_0 R) \le \theta_0^{2\alpha} \mathrm{Exc}(R).
\end{equation}
Indeed, if $r/R \ge \theta_0$, then the estimate~\eqref{excessstep1} is trivial at the expense of having possibly large constant on the right-hand side. Otherwise, let $n$ be the unique positive integer satisfying $\theta_0^{n-1} \le r/R < \theta_0^n$. Then, by induction,
\begin{equation}\nonumber
 \mathrm{Exc}(r) \lesssim \mathrm{Exc}(\theta_0^n R) \lesssim (\theta_0^{n})^{2\alpha} \mathrm{Exc}(R) = \theta_0^{2\alpha} (\theta_0^{n-1})^{2\alpha} \mathrm{Exc}(R) \lesssim  \left( \frac{r}{R} \right)^{2\alpha} \mathrm{Exc}(R),
\end{equation}
which concludes the proof of (\ref{excessstep1}) and thereby the proof of Theorem~\ref{prop1}.

\section{Proof of Lemma~\ref{lmcac}}

After possibly adding a constant to $u$ we may assume $c=0$.  Fix two radii $R>0$ and $0<\rho<\frac{R}{2}$, and suppose that $u$ is an $a$-harmonic function in $B_R$.  Let $\eta$ denote a smooth cutoff function such that $0\leq \eta \leq 1$ and
$$\eta(x)=\left\{\begin{array}{ll} 1 & \textrm{in}\;\;\overline{B}_{R-\rho}, \\ 0 & \textrm{in}\;\;\mathbb{R}^d\setminus B_R, \end{array}\right.$$
and satisfying
$$\abs{\nabla \eta}\lesssim \frac{1}{\rho}.$$
Using H\"older's inequality and \eqref{cac1}, we have
\begin{equation}\label{lmcac_11}\biggl(\fint_{B_{R-\rho}}\abs{\nabla u}^{\frac{2q}{q+1}}\biggr)^{\frac{q+1}{q}}\lesssim \Lambda\fint_{B_{R-\rho}} \nabla u\cdot a\nabla u \lesssim \Lambda  \fint_{B_R} \eta^2 \nabla u\cdot a\nabla u.\end{equation}
Then, by testing the equation $-\nabla\cdot a\nabla u$ against $\eta^2u$, and using the identity
$$\nabla(\eta^2 u)\cdot a\nabla u=\eta^2 \nabla u \cdot a\nabla u+2\eta u \nabla \eta\cdot a\nabla u,$$
it follows from the definition of $\mu$ in (\ref{f2}) that
$$\fint_{B_R} \eta^2\nabla u\cdot a\nabla u=-2\fint_{B_R}\eta u \nabla \eta \cdot a\nabla u\leq\fint_{B_R}2\eta \abs{a\nabla u}\abs{u\nabla\eta}\leq2\fint\mu^{\frac{1}{2}} (\eta^2 \nabla u\cdot a \nabla u)^{\frac{1}{2}}\abs{u\nabla \eta}.$$
Following an application of H\"older's inequality, we obtain
\begin{equation*}\fint_{B_R} \eta^2\nabla u\cdot a\nabla u \lesssim \biggl(\fint_{B_R}\eta^2\nabla u \cdot a\nabla u\biggr)^{\frac{1}{2}}\biggl(\fint_{B_R}\mu^p\biggr)^{\frac{1}{2p}}\biggl(\fint_{B_R}\abs{u\nabla \eta}^{\frac{2p}{p-1}}\biggr)^{\frac{p-1}{2p}}.\end{equation*}
Then, after dividing by the square-root of the left-hand side and using properties of the cutoff $\eta$, we have
\begin{equation}\label{lmcac_12}\fint_{B_R} \eta^2\nabla u\cdot a\nabla u \lesssim \frac{\Lambda}{\rho^2} \biggl(\fint_{B_R}\abs{u}^{\frac{2p}{p-1}}\biggr)^{\frac{p-1}{p}}.\end{equation}
In combination, inequalities (\ref{lmcac_11}) and (\ref{lmcac_12}) imply
$$\biggl(\fint_{B_{R-\rho}}\abs{\nabla u}^{\frac{2q}{q+1}}\biggr)^{\frac{q+1}{q}}\lesssim \Lambda\fint_{B_{R-\rho}} \nabla u\cdot a\nabla u \lesssim \frac{\Lambda^2}{\rho^2}\biggl(\fint_{B_R \setminus B_{R-\rho}}\abs{u}^{\frac{2p}{p-1}}\biggr)^{\frac{p}{p-1}},$$
which completes the proof of Lemma~\ref{lmcac}.

\section{Appendix}

\subsection{An Alternate Construction of the Flux Correction $\sigma$ in Lemma~\ref{lm1}.}

We claim that for any exponent in the range of
\begin{equation}\label{o2}
(2^*)'<r< 2,
\end{equation}
where $(2^*)'$ denotes the dual exponent of $2^*$, which in turn denotes the Sobolev exponent for $2$, we have the following existence result:
For a stationary random field $g$ with $\langle |g|^r\rangle^\frac{1}{r}<\infty$ there exists 
a curl-free stationary random field $\nabla\sigma$ of vanishing expectation $\langle\nabla\sigma\rangle=0$ with
\begin{equation}\label{o1}
-\triangle\sigma=\nabla\cdot g\quad\mbox{and}\quad\langle|\nabla\sigma|^r\rangle^\frac{1}{r}\lesssim\langle|g|^r\rangle^\frac{1}{r}.
\end{equation}
We note that the range (\ref{o2}) is sufficient for our purposes: We use it for $\frac{1}{r}=\frac{1}{2}+\frac{1}{2p}$
and note that under the weaker assumption $\frac{1}{d}>\frac{1}{2p}$, as compared with $\frac{1}{p}+\frac{1}{q}\le \frac{2}{d}$, we indeed
have (\ref{o2}) in the reciprocal form of $\frac{1}{(2^*)'}=1-(\frac{1}{2}-\frac{1}{d})
=\frac{1}{2}+\frac{1}{d}>\frac{1}{r}>\frac{1}{2}$.

\medskip

By applying the same reasoning as in the proof of Lemma \ref{lm1}, we may assume that $\langle|g|^2\rangle<\infty$, so that by
Riesz' representation theorem that there exists a curl-free stationary random field $\nabla\sigma$ with
$\langle\nabla\sigma\rangle=0$ and
\begin{equation}\nonumber
-\triangle\sigma=\nabla\cdot g\quad\mbox{and}\quad\langle|\nabla\sigma|^2\rangle^\frac{1}{2}\lesssim\langle|g|^2\rangle^\frac{1}{2}.
\end{equation}
By a standard duality argument, it is enough to establish (\ref{o1}) in the dual range of exponents, that is,
\begin{equation}\label{o6}
\langle|\nabla\sigma|^r\rangle^\frac{1}{r}\lesssim\langle|g|^r\rangle^\frac{1}{r}\quad\mbox{for}\quad
2< r<2^*.
\end{equation}
The first ingredient is that thanks to $r<2^*$ in conjunction with $\langle|\nabla\sigma|^2\rangle^\frac{1}{2}<\infty$
and $\langle\nabla\sigma\rangle=0$ and $r< 2^*$, a repetition of the proof leading Lemma \ref{lm2} we have the sublinearity
\begin{equation}\label{o5}
\lim_{R\uparrow\infty}\frac{1}{R}\langle\fint_{B_R}|\sigma-\fint_{B_R}\sigma|^r\rangle^\frac{1}{r}=0.
\end{equation}
The second ingredient is that thanks to $r> 2$ and $d\ge 2$, we have the following strengthening
of Calderon-Zygmund's estimate (in physical space)
\begin{equation}\label{o4}
\big(\int|\nabla\tilde\sigma|^r\big)^\frac{1}{r}\lesssim\big(\int|\tilde g|^r\big)^\frac{1}{r}+R\big(\int|\tilde f|^r\big)^\frac{1}{r},
\end{equation}
provided the functions $\tilde\sigma$, $\tilde f$ and the field $\tilde g$ satisfy
\begin{equation}\label{o3}
-\triangle\tilde\sigma=\nabla\cdot\tilde g+\tilde f\quad\mbox{and}\quad
{\rm supp}\,\tilde\sigma,\tilde g,\tilde f\subset B_R.
\end{equation}

\medskip

We consider $\tilde\sigma:=\eta(\sigma-c)$ for some constant $c$ and some smooth function $\eta$ supported
in $B_R$ to be fixed soon
and note that (\ref{o3}) is satisfied for 
\begin{equation}\nonumber
\tilde g:=\eta g+2(\sigma-c)\nabla\eta\quad\mbox{and}\quad \tilde f:=-\nabla\eta\cdot g+(\sigma-c)\triangle\eta.
\end{equation}
Choosing the cut-off of the form $\eta(x)=\hat\eta(\frac{x}{R})$ with $\hat\eta(\hat x)=1$ for $|\hat x|\le\frac{1}{2}$
we see that (\ref{o4}) turns into
\begin{equation}\nonumber
\big(\int_{B_{\frac{R}{2}}}|\nabla\sigma|^r\big)^\frac{1}{r}\lesssim\big(\int_{B_R}|g|^r\big)^\frac{1}{r}
+\frac{1}{R}\big(\int_{B_R}|\sigma-c|^r\big)^\frac{1}{r}.
\end{equation}
We now choose $c=\fint_{B_R}\sigma$ and apply $\langle(\cdot)^r\rangle^\frac{1}{r}$ to 
the above. By stationarity of $\nabla\sigma$ and $g$ we obtain
\begin{equation}\nonumber
\langle|\nabla\sigma|^r\rangle^\frac{1}{r}\lesssim\langle|g|^r\rangle^\frac{1}{r}
+\frac{1}{R}\big\langle\fint_{B_R}|\sigma-\fint_{B_R}\sigma|^r\big\rangle^\frac{1}{r}.
\end{equation}
With help of (\ref{o5}) we obtain (\ref{o6}).

\bibliographystyle{amsplain}
\bibliography{bella}

\end{document}